\documentclass[11pt]{article}
\usepackage[a4paper,margin=26mm]{geometry}
\usepackage{amsmath,amssymb,amsthm,booktabs,tabularx,microtype}
\usepackage{graphicx,float,xurl}
\usepackage{caption}
\usepackage[hidelinks,unicode]{hyperref}
\usepackage{enumitem}
\newtheorem{theorem}{Theorem}[section]
\newtheorem{lemma}[theorem]{Lemma}
\newtheorem{proposition}[theorem]{Proposition}
\newcommand{\repofile}[2]{\href{https://github.com/Jxxxx-2001/mordell-discriminant/blob/31cc934b7dfa0ef2308566a1bce48cbfb2e01559/#1}{#2}}
\newcommand{\repodir}[2]{\href{https://github.com/Jxxxx-2001/mordell-discriminant/tree/31cc934b7dfa0ef2308566a1bce48cbfb2e01559/#1}{#2}}
\providecommand{\doi}[1]{\href{https://doi.org/#1}{doi:\nolinkurl{#1}}}
\hypersetup{pdfauthor={Xiang Ji, Tiancheng Xu, Wensheng Yu}}
\author{Xiang Ji, Tiancheng Xu, Wensheng Yu\\[3pt]
  \small School of Electronic Engineering\\
  \small Beijing University of Posts and Telecommunications\\
  \small Beijing 100876, China}\date{}
\usepackage{titlesec}
\titleformat{\section}{\normalfont\Large\bfseries}{\thesection}{1em}{}
\titleformat{\subsection}{\normalfont\large\bfseries}{\thesubsection}{1em}{}
\titleformat{\subsubsection}{\normalfont\normalsize\bfseries}{\thesubsubsection}{1em}{}
\titlespacing*{\section}{0pt}{2.4ex plus .5ex minus .2ex}{1.2ex plus .2ex}
\titlespacing*{\subsection}{0pt}{2ex plus .4ex minus .2ex}{1ex plus .2ex}
\titlespacing*{\subsubsection}{0pt}{1.8ex plus .3ex minus .2ex}{.8ex plus .2ex}
\titlespacing*{\paragraph}{0pt}{2ex plus .3ex minus .2ex}{1em}
\titlespacing*{\subparagraph}{0pt}{2ex plus .3ex minus .2ex}{1em}

\title{Optimal Configurations for Mordell's Problem}
\begin{document}\maketitle
\begin{abstract}
Mordell's problem asks for the largest product of squared pairwise distances among $n$ points in the complex plane with a prescribed sum of squared distances from the origin. When this sum equals $n$, the maximum for $2\le n\le5$ is known to be $n^n$, attained precisely by the $n$th roots of unity up to rotation and relabeling. For $n\ge6$, Chalk showed that a regular polygon together with its center gives a larger product than the regular $n$-gon. We prove that for $6\le n\le9$ the maximum is $(n-1)^{n-1}r_n^{n(n-1)}$, where $r_n=\sqrt{n/(n-1)}$. Up to the same symmetries, the unique maximizer consists of the origin and the vertices of a regular $(n-1)$-gon of radius $r_n$. For six points, algebraic elimination reduces the proof to an exact comparison of stationary values. For seven points, we compare all twelve classes in a published classification of equilibria. For eight and nine points, rigorous interval arithmetic establishes an inequality for every admissible pair. Summing over pairs gives the sharp global bound, and the equality conditions determine the configuration. An explicit ten-point comparison and a known instability result together exclude Chalk's configuration for every $n\ge10$. For each $n$ from two through eight, Mordell's problem and the minimum-diameter problem studied by Bateman and Erd\H{o}s have a common optimal configuration.

\end{abstract}

\section{Introduction and main results}
Mordell's problem, in the sense used here, asks for the largest product of squared pairwise distances under a quadratic moment constraint. For $n\ge2$, write $N=\binom n2$ and set

\begin{equation}\label{objective} Q(z)=\sum_{i=1}^n|z_i|^2,\qquad \Delta(z)=\prod_{i<j}|z_i-z_j|^2,\qquad M_n=\max_{Q(z)=n}\Delta(z).\end{equation}
Let $P_n$ be the unit regular $n$-gon, with $P_2=\{-1,1\}$. For $n\ge3$ let $S_n$ consist of the origin and the vertices of a regular $(n-1)$-gon of radius $\sqrt{n/(n-1)}$ (an antipodal pair when $n=3$). Define $C_n=\Delta(S_n)$.

\begin{figure}[H]
\centering
\includegraphics[width=.77\linewidth]{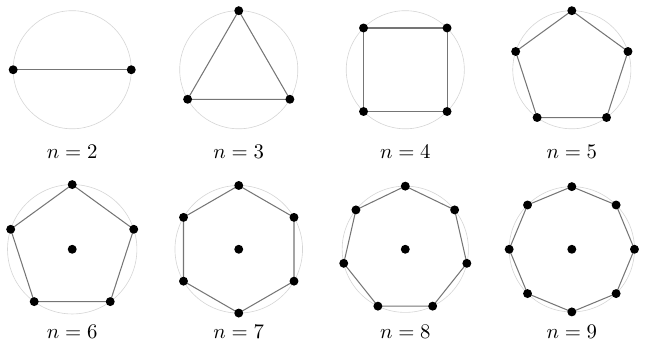}
\caption{Optimal configurations for Mordell's problem, $2\le n\le9$.}\label{fig:small}
\end{figure}

\begin{theorem}\label{main}
Under the normalization in \eqref{objective}, the sharp values and all maximizers up to rotation and relabeling are given in Table~\ref{tab:sharp}.

\end{theorem}
\begin{table}[H]
\centering
\caption{Sharp values and maximizing configurations for two through nine points.}\label{tab:sharp}
\begin{tabular}{ccc}\toprule
$n$ & Maximum $M_n$ & Maximizing configuration\\\midrule
$2,3,4,5$ & $n^n$ & $P_n$\\
$6,7,8,9$ & $(n-1)^{n-1}\left(n/(n-1)\right)^N$ & $S_n$\\\bottomrule
\end{tabular}
\end{table}
Figure~\ref{fig:small} shows the configurations in Theorem~\ref{main}. At six points a central point appears and the outer radius increases from $1$ to $\sqrt{6/5}$.

Birch's proposed optimizer was $P_n$. Mordell treated the three-point case and the restriction to a common unit circle \cite{mordell}. Chalk exhibited $S_n$ and showed that $C_n>n^n$ for every $n\ge6$ \cite{chalk}. This leaves a sharp question: when does Chalk's lower bound attain the maximum? Table~\ref{tab:sharp} answers it for $6\le n\le9$. Section~\ref{ten} explains its failure for every $n\ge10$.

The shapes have a history in vortex theory. Campbell and Ziff tabulated equilibria \cite{cz}. Gueron and Shafrir explicitly listed $S_6,\ldots,S_9$ among the minimizers suggested by numerical optimization \cite[Sec.~5]{gs}. Their constrained negative logarithmic energy at $w=z/\sqrt n$ equals $(N\log n-\log\Delta(z))/2$. Thus the numerical predictions concern exactly the shapes in Theorem~\ref{main}. For six and seven points, a publicly available manuscript also obtains the sharp inequalities using the classification of Faug\`ere and Svartz \cite{fs,jiang}. Below we distinguish inputs from that classification from comparisons actually certified here.

The proof has two parts. For six and seven points we reduce the maximization to a finite comparison of stationary values. For eight and nine points we instead construct a bound that applies to every possible maximizer at once. In both arguments the essential issue is coverage: all stationary branches in the first, and all admissible pairs in the second. The six-point argument uses a recorded exact quotient-dimension computation, whereas the seven-point argument uses the published classification in \cite{fs}. These dependencies are stated at the point of use.

Section 2 reviews the known cases and derives the stationary equations. Sections 3 and 4 prove the six- through nine-point results. Section~\ref{ten} excludes $S_n$ for every $n\ge10$. Section 6 compares the two extremal problems. Section~\ref{conclusion} summarizes how the maximizing configurations change with $n$ and how their optimality is proved. The appendices describe the computations. Their exact inputs and verification programs are archived in the fixed repository revision cited in Appendix~\ref{repository}.

Below, S6, S7, S8, S9, S10 and G identify the repository components mapped in Appendix~\ref{repository}; they are not files embedded in the TeX submission.

\section{Two through five points and common reductions}
\subsection{The known cases}
For two points, $\Delta=4-|z_1+z_2|^2$. For three, the squared distances sum to $3Q-|\sum z_i|^2\le9$. Their product is at most $27$ by the arithmetic-geometric mean inequality, with equality precisely for a centered equilateral triangle. For four points we invoke Lin's theorem \cite{lin}; a short alternate proof changes Hadamard coordinates so the four-point Vandermonde becomes the three-point Vandermonde of three squared coordinates \cite[Appendix~C]{wang}. For five, Wang reduces the stationary equations to polynomial branches, establishes reflection symmetry and compares the five physical shapes \cite{wang}. Xu, Dou and Yu use Wang's coefficient conditions but directly bound the discriminant on the two generic quartic branches \cite{xu}. These low-degree theorems give the first line of Theorem~\ref{main}. We use their equality statements as well as their upper bounds.

\subsection{Two attainable values and the stationary equations}
For a monic polynomial $p$ with distinct roots, $\Delta=|\prod_i p'(z_i)|$. Applied to $X^n-1$ this gives $\Delta(P_n)=n^n$. At the roots of $X^n-X$ the absolute derivatives are $1$ once and $n-1$ at each of the other roots. Scaling the roots gives the value $C_n$ in Table~\ref{tab:sharp}.

In particular, writing $N=\binom n2$,
\begin{equation}\label{candidatevalues}
\Delta(P_n)=n^n,\qquad C_n=(n-1)^{n-1}\left(\frac n{n-1}\right)^N.
\end{equation}
These computations supply lower bounds. To prove optimality, we must show that an arbitrary maximizer cannot do better. The following reduction imposes no symmetry assumption.

\begin{lemma}\label{stat}
A maximizer exists, is centered, has distinct points, and satisfies

\begin{equation}\label{forces}\sum_{j\ne i}\frac{1}{z_i-z_j}=\lambda\bar z_i,\qquad \lambda=(n-1)/2.\end{equation}
\end{lemma}
\begin{proof}
The sphere $Q=n$ is compact and $\Delta(P_n)>0$. Centering preserves all distances and changes $Q$ to $q=n-n|\sum_i z_i/n|^2$. If $q<n$, rescaling multiplies the positive product by $(n/q)^N>1$. The real Lagrange multiplier equations for $\log\Delta$ give \eqref{forces} with a real multiplier. Multiply by $z_i$ and sum: each pair contributes one, so $\lambda n=N$.

\end{proof}
\section{Six and seven points: finite stationary comparisons}
\subsection{Six points: exhausting the stationary branches}
Our aim is to show that every maximizing configuration other than $S_6$ has value at most $6^6<C_6$. It is enough to derive a finite list containing all collision-free stationary configurations and to bound every branch on that list. Some algebraic solutions may fail to represent complex points satisfying the conjugacy conditions; retaining such solutions is harmless when their objective bound is already smaller.

\paragraph{Necessary coefficient equations.}
Center and rescale the six roots to $w_i$. Write $q=\sum|w_i|^2$, $\kappa=30/q$ and $p=X^6+AX^4+bX^3+cX^2+dX+e$. The quantities $Y_i=p''(w_i)/p'(w_i)=\kappa\bar w_i$ are distinct. Conjugating the force equations gives $2\sum_{j\ne i}(Y_i-Y_j)^{-1}=w_i$. To eliminate the pairing of $w_i$ and $Y_i$, put

\[
R(Y,T)=\operatorname{Res}_X(p,(Y+TX)p'-p''),\quad R_0=R(Y,0),\quad D=[Y^6]R_0.
\]
Writing $R_T=\partial R/\partial T$, the product formula gives $R_0=D\prod(Y-Y_i)$ and $R_T(Y_i,0)=w_iR_0'(Y_i)=R_0''(Y_i)$. Since the difference has degree below six, $R_T(Y,0)=R_0''(Y)$. Comparing coefficients gives four necessary polynomial equations. At a physical solution the coefficients of $R_0/D$ are the conjugate coefficients of $p$ multiplied by the corresponding powers of $\kappa$. This relation also covers zero coefficients.

The coefficient cases are exhaustive. If $A\ne0$, scale to $A=1$ and first separate $b\ne0$ from $b=0$; in the latter case separate $d\ne0$ from $d=0$. If $A=0$, use the first nonzero coefficient among $b,c,d$, or else obtain $X^6+e$. Scaling a nonzero coefficient to one is a complex similarity and preserves the normalized objective. Thus no reflection or collinearity assumption enters the branch division.

\paragraph{Completeness of the generic branch.}
If $Ab\ne0$, set $A=1$, $B=b^2$ and $k=d/b$. After division by the permitted nonzero factors, the four coefficient equations are polynomial in $B,c,k,e$. Adjoin $tBD-1=0$ to exclude $B=0$ and collisions. Let $I$ be the ideal they generate in $\mathcal R=\mathbb Q[t,e,k,c,B]$. Repository material S6 specifies a triangular list consisting of six square-free, pairwise coprime factors in $B$, of degrees $1,1,2,6,11,14$, and polynomial expressions for the other coordinates.

\begin{lemma}\label{sixcomplete}
The triangular list contains every complex solution of the generic necessary equations.

\end{lemma}
\begin{proof}
Let $J$ be the intersection of the six graph ideals defined by the list. Exact substitution in the six number fields verifies every generator of $I$, including $tBD-1$, so $I\subseteq J$. The coprime, square-free factors give $\dim_{\mathbb Q}\mathcal R/J=35$. An independent rational Magma computation gives $\dim_{\mathbb Q}\mathcal R/I=35$; its exact input and recorded output are available in S6. The natural surjection $\mathcal R/I\to\mathcal R/J$ is therefore an isomorphism. Hence $I=J$.

\end{proof}
The dimension computation is a dependency of this proof. The repository routine rechecks the field identities and subsequent comparisons, but reads the saved dimension result; the separate Magma input permits that computation to be repeated. Merely checking that the listed solutions satisfy the equations would not establish the lemma.

\paragraph{Conjugacy and objective comparison.}
At a physical solution, $a=[Y^4]R_0/D=\kappa^2$ is real and positive. Since $\sum w_i^2=-2$, we have $q\ge2$ and $a\le225$. Equality forces every $w_i$ to be imaginary, hence $B=b^2\le0$. The branch with $a=225$ instead has $B=32/45$, and is excluded.

Reality of $a$ requires care: it does not automatically imply that all coefficients are real. In the first five factors, exact field computations give $\mathbb Q(B)=\mathbb Q(a)$, so real $a$ does force real $B,c,k,e,D$. The last factor has degree fourteen, while $a$ has degree seven. There a computed involution $\sigma$ fixes $a$ and must act as complex conjugation at a physical solution. Its remaining admissible real values of $a$ make $B$ and $\sigma(B)$ distinct real numbers. They therefore cannot be conjugates. This eliminates that factor without assuming reflection symmetry. The identities and isolating intervals are given in Appendix~\ref{algebra} and S6.

Four real roots remain, one each from factors 0 and 2 and two from factor 4 in Appendix~\ref{algebra}. If $\Phi_6=(6/q)^{15}\Delta(w)$ is the normalized objective, then

\[
\Phi_6^2=\frac{a^{15}D^2}{5^{30}}.
\]
Rational interval evaluation gives a strict upper bound below $6^{12}$ at each of these roots. Table~\ref{tab:bounds} gives the four bounds. Thus none of the generic branches can maximize the objective.

\paragraph{Zero-coefficient cases.}
Table~\ref{tab:exceptional} lists the exceptional cases after scaling the first nonzero coefficient to one. These are rational ideal computations in S6; the collision factor is saturated throughout. A unit ideal means that the chart has no collision-free solution.

\begin{table}[H]
\centering
\caption{Exceptional coefficient cases for six points.}\label{tab:exceptional}
\small\begin{tabular}{ll}\toprule
Chart & Consequence\\\midrule
$A=1,\ b=0,\ d\ne0$ & unit ideal\\
$A=0,\ b=1$ & $c=d=0,\quad25e^2-17e+1=0$\\
$A=b=0,\ c=1$ & unit ideal\\
$A=b=c=0,\ d=1$ & $e=0$\\
$A=b=c=d=0$ & $X^6+e$\\
$A=1,\ b=d=0$ & $X^6+X^4+uX^2+v$\\\bottomrule
\end{tabular}
\end{table}
For the even chart the eliminant is $(5u-1)(5u^2+36u-9)(5u^3-9u^2+15u-3)$. The conjugacy constraints force real $u,v$; all four real roots give a squared objective below $6^{12}$. In the $b=1$ chart the roots form two equilateral triples, with polynomial $(X^3-U)(X^3-V)$. After rescaling back to $Q=6$, write $|U|=r^3,|V|=s^3$, so $r^2+s^2=2$. With $x=rs\in[0,1]$, one has $[rs(r^3+s^3)]^2=2x^2(2-x)^2(1+x)\le4$ by differentiation. Thus $\Delta=3^6|UV|^2|U-V|^6\le6^6$. The remaining shapes are $P_6$ and $S_6$, with $6^6<C_6$. The branch division above covers every stationary maximizer, so this proves the six-point assertion, including uniqueness.

\begin{theorem}\label{six}
Under $Q=6$, one has $\Delta\le C_6$, with equality exactly for $S_6$ up to rotation and relabeling.

\end{theorem}
\subsection{Seven points: comparison using a complete classification}
Here a published classification supplies the exhaustive list. The task remaining is to identify enough genuine equilibria to fill that list and to order their objective values rigorously.

Set $u_i=\sqrt3\,z_i$. Equation \eqref{forces} becomes $\sum_{j\ne i}(u_i-u_j)^{-1}=\bar u_i$. We use the full physical classification of \cite[Sec.~5.7 and Fig.~4]{fs}: there are twelve classes up to rotation, reflection and permutation. We use the number of classes only. It suffices to certify twelve genuine configurations with pairwise disjoint value intervals. No matching against a drawing is required.

In real coordinates let $F_i=\sum_{j\ne i}(z_i-z_j)/|z_i-z_j|^2-3z_i$. Set $y_0=0$, $x_0>0$ and omit the $y_0$ equation. The identity $\sum_i(x_iF_{iy}-y_iF_{ix})=0$ recovers it. Summing the forces and taking their scalar product with the positions recover centering and $Q=7$. On each of ten rational boxes of radius $r=2^{-80}$, S7 checks, for midpoint $m$ and preconditioner $C$,
\[
\|Cf(m)\|_\infty<r/4,\qquad \sup_X\|I-CDf\|_\infty<1/4,
\]
where $f$ denotes the retained equations. Thus $x\mapsto x-Cf(x)$ maps the closed box into itself and is a contraction. Also $C$ is invertible, since $CDf$ is within norm $1/4$ of the identity. Its fixed point is therefore a genuine solution of $f=0$. Every box is separated from collisions and satisfies $x_0>0$. This proves existence, rather than just a small residual at a numerical approximation.

Interval multiplication of the 21 squared distances gives ten pairwise disjoint objective intervals, all strictly below $7^7$. The intervals are displayed in Appendix~\ref{algebra}. Distinct objective values imply distinct classes even after allowing reflection.

The other two orbits are explicit stationary configurations: $P_7$ gives $823543$, and $S_7$ gives $7^{21}/6^{15}>823543$. Hence only $S_7$ attains the maximum. Its reflected image lies in the same rotation and permutation orbit.

We have thus certified twelve different classes: ten interval solutions and the two explicit configurations. The external twelve-class theorem leaves no additional class in which a maximizer could lie.

\begin{theorem}\label{seven}
Under $Q=7$, one has $\Delta\le C_7$, with equality exactly for $S_7$ up to rotation and relabeling. This conclusion uses the complete equilibrium classification of Faug\`ere and Svartz.

\end{theorem}
\section{Eight and nine points: global bounds from pair inequalities}
Here we bound all pairs of any maximizing configuration. The force equations enter only through necessary conditions. In particular, a point at the center and equally spaced outer angles will be consequences of equality.

The proof can be read in three steps. First we find a pairwise upper bound whose sum is at most the value of the centered polygon. Next we prove that every pair at an arbitrary global maximizer belongs to a domain on which that bound is valid. Finally, equality in the summation determines the radii and angles. The domain may contain triples that no actual configuration realizes; this only makes the inequality stronger than required.

\subsection{A common summation argument}
To make the candidate outer radius equal to one, put $m=n-1$, $x_i=\sqrt{m/n}\,z_i$, $u_i=|x_i|^2$, $t_{ij}=\operatorname{Re}(x_i\bar x_j)$ and $d(u,v,t)=u+v-2t$. Thus $\sum u_i=m$. Use $s=4$ for $n=8$ and $s=6$ for $n=9$, and the column vectors

\[
 F(u)=(u-m/n,u(u-1),\ldots,u^{s-1}(u-1))^T,\quad
 V(u)=(1,u-1,\ldots,(u-1)^{s-1})^T .
\]

The first vector records radial constraints, while the second will record angular moments. The factor $u-1$ makes these constraints vanish on the unit outer ring. Their purpose is to turn the sum of a two-point kernel into a nonnegative quadratic form in moments of the whole configuration.

Let $H_0=1,H_1=t,H_k=2tH_{k-1}-uvH_{k-2}$. For $u=|x|^2,v=|y|^2,t=\operatorname{Re}(x\bar y)$, the recurrence gives $H_k=\operatorname{Re}((x\bar y)^k)$. For real symmetric matrices $Q_0,\ldots,Q_{m-1}$ of order $s$, define

\begin{align}
 K(u,v,t)&=F(u)^TQ_0F(v)+\sum_{k=1}^{m-1}H_k V(u)^TQ_k V(v),\nonumber\\
 p(u)&=\big(nu/2+b-K(u,u,u)/2\big)/m,\quad b=(m\log m-N)/n,\nonumber\\
 h(u,v,t)&=p(u)+p(v)-K(u,v,t).\label{kernel}
\end{align}

The diagonal correction in $p$ is chosen so that all single-point terms cancel when $h$ is summed. Consequently the following lemma reduces the global extremal problem to one inequality in the three real variables $u,v,t$.

\begin{lemma}\label{sum}
Suppose $Q_k\succ0$ for all $k$ and $h(u_i,u_j,t_{ij})\ge\log |x_i-x_j|^2$ for every distinct pair at every global maximizer. Then $M_n=C_n$, with equality precisely on the orbit of $S_n$.

\end{lemma}
\begin{proof}
Write $K_{ij}=K(u_i,u_j,t_{ij})$, $L_0=\sum_i F(u_i)$ and $L_k=\sum_i x_i^k V(u_i)$. Pairing conjugates and summing \eqref{kernel} give

\begin{equation}\label{sumid}
\sum_{i,j}K_{ij}=\sum_{k=0}^{m-1} L_k^*Q_kL_k\ge0,\quad
\sum_{i<j}h(u_i,u_j,t_{ij})=m\log m-\tfrac12\sum_{i,j}K_{ij}.
\end{equation}
Thus $\log\Delta(x)\le m\log m$ at every global maximizer. Since a maximizer exists, this bounds the original maximum by $C_n$ after rescaling; $S_n$ already attains that value. Equality forces every $L_k$ to vanish. Subtracting the second component of $L_0$ from the third gives $\sum u_i(u_i-1)^2=0$. Hence exactly $m$ radii squared are one and the remaining radius is zero. For $1\le k<m$, the first component of $L_k$ gives $\sum x_i^k=0$. Newton's identities make the monic polynomial of the $m$ outer roots $X^m+c$ with $|c|=1$. This proves the equality assertion.

\end{proof}
\subsection{Necessary conditions at any global maximizer}
We now obtain a bounded domain containing every pair that the lemma needs. The second variation supplies separation and moment bounds. These are necessary conditions even when the maximum is degenerate.

Set $A_{ij}=-(z_i-z_j)^{-2}$ for $i\ne j$ and choose the diagonal for zero row sums. Then $Az=\lambda\bar z$. The second variation of $\log\Delta-\lambda Q$ is $-2\operatorname{Re}(v^TAv)-2\lambda\|v\|^2$. It is nonpositive on the tangent space, and the radial vector $z$ is orthogonal to this space for the associated real quadratic form, with negative radial value. Hence it is nonpositive everywhere. Applying the estimate at $v$ and $iv$ bounds the spectrum of the real symmetric operator $v\mapsto\overline{Av}$ by $[-\lambda,\lambda]$. In particular $\|A\|\le\lambda$. Centering the coordinate vectors before testing entries gives $|A_{ij}|\le \lambda(1-1/n)$, so $d\ge2/m$ in the $x$ scale. None of these conclusions requires a strictly definite Hessian.

\subsection{Eight points: a radial bound and its certificate}\label{eightdomain}
Write $R_j=\sum|z_i|^j$. Acting on the centered coordinate powers of orders two and three gives $(Az^2)_i=7(|z_i|^2-1)$ and $(Az^3)_i=(21/2)z_i|z_i|^2-13z_i$. The norm bound yields $R_4\le32/3$ and $98R_6-273R_4+1352\le0$ (a nonnegative term has been dropped). For $s_0=\max|z_i|^2$, Cauchy's inequality gives $R_6\ge s_0^3+(R_4-s_0^2)^2/(8-s_0)$. The resulting left side decreases with $R_4$ up to $32/3$: its stationary value in the $R_4$ variable occurs at $s_0^2+39(8-s_0)/28>32/3$ for $s_0\ge1$. Consequently

\[
882s_0^3-2352s_0^2+1755s_0-1496\le0.
\]
Sturm isolation in S8 gives just one root on $[1,8]$, smaller than $8U_8/7$, where $U_8=965267/524288$. Thus $0\le u,v\le U_8$, $t^2\le uv$ and $d\ge2/7$ include every required pair.

\begin{proposition}\label{cert8}
The matrices specified in S8 are strictly positive definite, and their function $h$ satisfies $h(u,v,t)\ge\log(u+v-2t)$ throughout

\[
\mathcal D_8=\{(u,v,t):0\le u,v\le U_8,\ t^2\le uv,\ u+v-2t\ge2/7\}.
\]
\end{proposition}
This is the computer-assisted step. The exact matrices, local contact tests and continuous-domain verification are specified in Appendix~\ref{repro}; the verification starts from the full box containing $\mathcal D_8$. Together with the necessary-domain argument above, Proposition~\ref{cert8} supplies every hypothesis of Lemma~\ref{sum}.

\begin{theorem}\label{eight}
Under $Q=8$, one has $\Delta\le C_8$, with equality exactly for $S_8$ up to rotation and relabeling.

\end{theorem}

\subsection{Nine points: additional radial moment constraints}
For nine points the certificate uses more information than a separate bound on each radius. The extra conditions relate the two radii of a tested pair to the seven radii that remain. We derive them before stating the inequality to be verified.

Write $T_j=\sum u_i^j$. The projected-power norm inequalities, now through order four, yield

\begin{equation}\label{mom}
T_2\le256/27,\quad 2T_3-5T_2+200/9\le0,\quad405T_4-1008T_3+912T_2-3072\le0 .
\end{equation}
For the last inequality, before rescaling one has $(Az^4)_i=16z_i^2|z_i|^2-21z_i^2-\sum z_j^2$. In the squared norm use $4\sum |z_j|^2z_j^2=(15/2)\sum z_j^2$ from the forces; the second-power test gives $|\sum z_j^2|^2\le324-27R_4$. Substitution gives $240R_8-672R_6+684R_4-2916\le0$, which rescales to the last part of \eqref{mom}.

Put $\rho=(225-\sqrt{18465})/90$, $L=(14372-100\sqrt{18465})/405$ and $f(a)=a^2(a-\rho)^2$. Add $504-405\rho>0$ times the middle inequality in \eqref{mom} to the last. The result is $\sum f(u_i)\le L$. In particular each $u_i\le(\rho+\sqrt{\rho^2+4\sqrt L})/2<U_9=1860281/1048576$. Besides $t^2\le uv$ and $d\ge1/4$, the verifier uses two necessary inequalities. With $a=u^2+v^2$, $s_0=u+v$, set

\[
W=\min\{256/27,\ a+5(8-s_0)/4\}.
\]
Cauchy's inequality for the other seven squared radii bounds their third moment below by $(T_2-a)^2/(8-s_0)$. Minimizing the resulting quadratic over $T_2\le256/27$ gives the first test below. For the second, $f(y)\ge(\rho^2/4)\operatorname{dist}(y,\{0,\rho\})^2$ for $y\ge0$. Round each of the other seven radii squared to its nearest endpoint $0$ or $\rho$ and use Cauchy's inequality for the seven errors. Thus every actual pair satisfies

\begin{align}
2[u^3+v^3+(W-a)^2/(8-s_0)]-5W+200/9&\le0,\nonumber\\
 f(u)+f(v)+(\rho^2/28)\min_{0\le k\le7,\ k\in\mathbb Z}(8-u-v-k\rho)^2&\le L.\label{pair9}
\end{align}

\begin{proposition}\label{cert9}
Let $\mathcal D_9$ consist of the triples with $0\le u,v\le U_9$, $t^2\le uv$, $u+v-2t\ge1/4$, and both inequalities in \eqref{pair9}. The matrices specified in S9 are strictly positive definite, and $h(u,v,t)\ge\log(u+v-2t)$ on $\mathcal D_9$.

\end{proposition}
The certificate also treats the boundary $uv=t^2$, including the antipodal contact of the outer octagon. At these contacts an ambient Hessian test is replaced by a one-sided argument on the physical domain; Appendix~\ref{repro} gives the correction explicitly. Thus boundary pairs are covered as well as interior pairs.

Every pair at a global maximizer belongs to $\mathcal D_9$ by the preceding derivation. Proposition~\ref{cert9} and Lemma~\ref{sum} now give the final small-cardinality case.

\begin{theorem}\label{nine}
Under $Q=9$, one has $\Delta\le C_9$, with equality exactly for $S_9$ up to rotation and relabeling.

\end{theorem}

The equality proof is the same in both cases. Strict positivity of the kernel matrices first fixes the radii, and the vanishing angular moments then fixes the polygon. No central point, reflection symmetry or equal angular spacing was assumed in the necessary-domain analysis.

\section{Failure of Chalk's configuration for ten or more points}\label{ten}
\subsection{Ten points: a strict comparison}
The preceding results establish global optimality of $S_n$ for four consecutive cardinalities. At ten points, an explicit feasible configuration already gives a larger product. This counterexample and its two-rectangle geometry occur in Gueron and Shafrir's 1999 paper \cite[Fig.~1 and pp.~13--14]{gs}. Here we verify the strict comparison directly from their printed coordinates.

\begin{proposition}\label{tenresult}
One has $M_{10}>C_{10}$.

\end{proposition}
\begin{proof}
Treat the printed decimal coordinates as exact rationals and clear their common denominator. This gives

\[
W=\{(\pm11183,0),\;(\pm32765,\pm14669),\;(\pm12323,\pm31572)\},
\]
with independent signs for each rectangle. Set $q=\sum_{w\in W}|w|^2=9999608174$ and $z=\sqrt{10/q}\,W$. Then $Q(z)=10$. Enumerating the 45 squared distances gives the integer inequality

\begin{equation}\label{historical}
9^{36}\prod_{i<j}|w_i-w_j|^2>q^{45}.
\end{equation}
It follows that $\Delta(z)>10^{45}/9^{36}=C_{10}$. Repository material S10 verifies \eqref{historical} using integer arithmetic only.

\end{proof}

This comparison uses only the feasibility of the normalized coordinates. It gives a strict lower bound for $M_{10}$ and leaves the maximizing configuration undetermined.

\subsection{Eleven or more points: instability}
For $Q(z)=n$ and $w=z/\sqrt n$, Gueron and Shafrir's constrained energy on the unit sphere is

\[
F(w)=-\sum_{i<j}\log|w_i-w_j|
     =\tfrac12\bigl(N\log n-\log\Delta(z)\bigr).
\]
For every $n\ge11$, Theorem~1(2) of \cite{gs} gives an unstable direction at the polygon with its center. On the constraint sphere there are therefore arbitrarily close configurations with strictly smaller energy. Scaling back to $Q=n$ gives a configuration with $\Delta>C_n$.

Together with Proposition~\ref{tenresult}, this proves $M_n>C_n$ for every $n\ge10$. Among the cases $n\ge6$ in Chalk's comparison, $S_n$ is thus globally optimal precisely for $6\le n\le9$.

\section{An Erd\H{o}s problem on minimum diameter}
Bateman and Erd\H{o}s \cite[pp.~306--307]{be} asked for the least diameter of a set of $n$ planar points whose mutual distances are at least one. For any configuration of distinct points, scaling by the reciprocal of its shortest distance makes that distance one and its diameter equal to

\[
 R(z)=\frac{\max_{i<j}|z_i-z_j|}{\min_{i<j}|z_i-z_j|}.
\]
Thus the comparison is between shapes up to similarity. Bateman and Erd\H{o}s give the values through six points and prove the seven-point case in their Theorem~3 \cite{be}. Bezdek and Fodor determine the eight-point minimum \cite{bf}. The shapes $P_2,\ldots,P_5,S_6,S_7,S_8$ attain those minima as well as the Mordell maxima in Theorem~\ref{main}. Their squared ratios, computed from the chord lengths, are respectively

\[
1,\quad 1,\quad 2,\quad (3+\sqrt5)/2,\quad(5+\sqrt5)/2,\quad4,\quad 1/(2-2\cos(\pi/7)).
\]
\begin{proposition}\label{ratio}
At eight points the two sets of optimizing shapes differ. Nine is the first cardinality with no common optimizer. For every $n\ge10$ both objectives exclude $S_n$.

\end{proposition}
\begin{proof}
Move the center of the unit regular heptagon a distance $1/10$. All its distances to the vertices are in $[0.9,1.1]$. The shortest outer chord is $2\sin(\pi/7)<0.9$ and the longest is $2\cos(\pi/14)>1.1$. Thus $R$ stays optimal, yet the displaced point cannot be the center of the circumcircle of the seven hull vertices. By Theorem~\ref{main} the resulting shape fails to maximize $\Delta$.

For nine points, a single explicit competitor is enough; there is no need to know the optimizer of the ratio problem. Use the integer witnesses in Appendix~\ref{integers}. At nine points $R(S_9)^2=4+2\sqrt2>925/136$. Since $S_9$ is the unique maximizing shape for $\Delta$, there is no common optimizer. The ten- and eleven-point witnesses give squared ratios $733/89$ and $289/29$, respectively smaller than $R(S_{10})^2=1/[4\sin^2(\pi/18)]$ and $R(S_{11})^2=\csc^2(\pi/10)$. Only attainable comparisons are needed.

The uniform construction for $n\ge12$ is given in Appendix~\ref{integers}. Together with Section~\ref{ten}, it proves the last assertion.

\end{proof}
The two transitions are thus distinct. For $n\ge10$ the exclusion of one shape by both criteria gives no identification of their actual optimizing sets.

\section{Conclusion}\label{conclusion}
The small-cardinality cases exhibit two successive maximizing shapes. The known two- through five-point results select the unit regular polygon, with an antipodal pair for two points; the six- through nine-point results select a regular polygon together with its center. The proofs for six and seven points compare stationary values, with completeness supplied by exact elimination in the former and a published classification in the latter. For eight and nine points, a pair inequality on a necessary domain gives a global bound without classifying every stationary configuration. Its equality conditions recover both the radii and the angular spacing.

For every $n\ge10$, the strict inequality $M_n>C_n$ excludes Chalk's configuration, by a direct comparison at ten points and by instability thereafter. Determining the sharp values and all maximizing configurations in this range remains open in the present treatment. The comparison with the minimum-diameter problem gives a different threshold. The two problems have a common optimal shape through eight points, but none for nine.

The code and verification materials are available at \url{https://github.com/Jxxxx-2001/mordell-discriminant}. Appendix~\ref{repository} identifies the fixed revision and the verification entry points.

\medskip
\noindent\textbf{Acknowledgments.} GPT-6 Astra assisted with calculations, proof development, and editorial work. The computational checks use exact algebra and rigorous interval arithmetic, with their external dependencies stated explicitly.

\clearpage
\appendix
\section{Algebraic branches and stationary-value comparisons}\label{algebra}
\subsection{The six-point equations and conjugacy test}
With the resultant notation of Section 3.1, put $r_j=[Y^j]R_0$ and $R_T=\partial_T R(Y,0)$. The identities $[Y^5]R_0=[Y^5]R_T=0$ and $[Y^4]R_T=30D$ leave the following four equations.

\begin{align*}
E_0&=[Y^3]R_T=0,& E_1&=[Y^2]R_T-12r_4=0,\\
E_2&=[Y]R_T-6r_3=0,& E_3&=[1]R_T-2r_2=0.
\end{align*}
Their expanded coefficients and the saturated quotient system are available in S6. Table~\ref{tab:factors} lists the field degree of $a=r_4/D$, the total number of its real roots, and the number remaining after the physical necessary tests. The factor index is the zero-based index in the data.

\begin{table}[H]
\centering
\caption{Real roots retained from the six generic factors for six points.}\label{tab:factors}
\begin{tabular}{rrrrr}\toprule
Factor & $\deg B$ & $\deg a$ & Real $a$ & Retained\\\midrule
0&1&1&1&1\\
1&1&1&1&0\\
2&2&2&2&1\\
3&6&6&2&0\\
4&11&11&3&2\\
5&14&7&3&0\\\bottomrule
\end{tabular}
\end{table}
For the last factor, the involution is determined by

\[
\sigma(B)=\frac{r_3^2D}{r_4^3},\quad
\sigma(c)=\frac{D r_2}{r_4^2},\quad
\sigma(k)=\frac{D r_1}{r_4r_3},\quad
\sigma(e)=\frac{D^2r_0}{r_4^3}.
\]
The formulas are interpreted in the corresponding number field. The checker verifies the inverses, $\sigma^2=1$, and $\sigma(a)=a$. The coefficient conjugacy relation implies that $\sigma$ must be conjugation at a physical solution. In the quadratic extension over $\mathbb Q(a)$, both $B+\sigma(B)$ and $B\sigma(B)$ belong to $\mathbb Q(a)$. At every remaining admissible real root, exact interval evaluation gives $(B-\sigma(B))^2>0$. Thus the two roots of that real quadratic are distinct and real, contradicting conjugacy.

For the four retained roots, the integer bounds in Table~\ref{tab:bounds} certify $\Phi_6^2<6^{12}$. The root intervals, stored as exact rational endpoints, are refined to width at most $10^{-70}$.

\begin{table}[H]
\centering
\caption{Certified upper bounds for $\Phi_6^2$ at the four retained roots.}\label{tab:bounds}
\begin{tabular}{rr}\toprule
Factor & Strict upper bound for $\Phi_6^2$\\\midrule
0&350000000\\
2&1803000000\\
4&224000000\\
4&10100000\\\bottomrule
\end{tabular}
\end{table}

\subsection{Seven-point value intervals}
Each row of Table~\ref{tab:seven} gives an open interval containing the objective of one certified equilibrium. Decimal endpoints denote exact rationals. The exact coordinate boxes have radius $2^{-80}$ and are available in S7.

\begin{table}[H]
\centering
\caption{Certified discriminant intervals for ten seven-point equilibria.}\label{tab:seven}
\begin{tabular}{rrr}\toprule
Class & Lower endpoint & Upper endpoint\\\midrule
1&151.332&151.333\\
2&9263.849&9263.850\\
3&40252.802&40252.803\\
4&40253.389&40253.390\\
5&118016.851&118016.852\\
6&159665.670&159665.671\\
7&268961.391&268961.392\\
8&486126.791&486126.792\\
9&666965.240&666965.241\\
10&822628.729&822628.730\\\bottomrule
\end{tabular}
\end{table}
The two explicit values are $\Delta(P_7)=823543$ and $\Delta(S_7)=7^{21}/6^{15}$. In particular, Table~\ref{tab:seven} separates all ten certified classes from each other and from both explicit ones.

\section{Kernels and verification of the continuous domains}\label{repro}
S8 and S9 specify the upper triangles of the matrices in \eqref{kernel} by exact expressions. In S8 the scalar basis is $1$ and the three numbers $\log(2-2\cos(2\pi j/7))$, $1\le j\le3$, over $\mathbb Q(\cos(2\pi/7))$. In S9 it is $1,\log(2-\sqrt2),\log2$ over $\mathbb Q(\sqrt2)$. The checkers reconstruct the polynomial $h$ coefficient by coefficient (193 and 444 monomials, respectively). Interval $LDL^T$ factorization certifies strict positivity of all matrix pivots. The archived rational and algebraic coefficients fix the kernels independently of any optimization used to find them.

Let $g=h-\log d$. The contact triples are $(0,1,0)$, its exchange, and $(1,1,\cos(2\pi j/m))$, $1\le j\le\lfloor m/2\rfloor$. At ordinary interior contacts exact value and gradient vanish, and interval Hessian bounds certify local nonnegativity. For S8 the center contact has $g_u>0$, $g_v=g_t=0$, so the same Taylor argument works for $u\ge0$. In S9, at the center and antipodal contacts, put $q_0=uv-t^2$. The checker subtracts $\tau q_0-q_0^2$ with $\tau-q_0>0$ on the box. The corrected residual has zero value and gradient and positive definite Hessian there; adding back $q_0(\tau-q_0)\ge0$ proves the required one-sided assertion.

The initial box is $[0,U_n]^2\times[-U_n,U_n]$. By symmetry of $u,v$ it suffices to cover $u\le v$; a box is discarded for this reason only when its lower $u$ endpoint exceeds its upper $v$ endpoint. In particular, the diagonal is retained. All other boxes must be proved outside the necessary domain, contained in a certified contact neighborhood, or positive by a Bernstein coefficient test; failing these tests the box is bisected. For $d_0>0$, a tangent bounds $\log d$ above. S9 can also use the cubic $\log d_0+y-y^2/2+y^3/3$ for $y=d/d_0-1>-1$ (negative fourth derivative). Subtracting such an upper polynomial from $h$, a positive lower bound for every tensor Bernstein coefficient proves positivity throughout the box. S9 permits the same one-sided corrections using $q_0$ with separately checked signs. Subdivision encloses every coefficient with outward rounding and checks the bisection midpoint exactly.

\begin{table}[H]
\centering
\caption{Box counts for the eight- and nine-point interval certificates.}\label{tab:coverage}
\begin{tabular}{rrrrr}\toprule
$n$ & visited & splits & terminal & unresolved\\\midrule
8&386915&193457&193458&0\\
9&31613&15806&15807&0\\\bottomrule
\end{tabular}
\end{table}
Table~\ref{tab:coverage} records the counts from the two complete runs. Every terminal action has a mathematical certificate. Subdivision uses closed child boxes with a common boundary. Thus no points are lost at a split, and an empty frontier covers the full initial box after the stated domain restrictions. A time limit or a nonempty frontier is a failure to certify, never a successful conclusion. These checks establish Propositions~\ref{cert8} and~\ref{cert9}; the analytic argument in Section 4 then gives the global bounds and all equality cases. This is a computer-assisted proof using exact algebra, Arb inclusions and outward-rounded IEEE binary64 operations.

\subsection{Repository inputs and verification requirements}\label{repository}
All material identifiers refer to the fixed repository revision \href{https://github.com/Jxxxx-2001/mordell-discriminant/tree/31cc934b7dfa0ef2308566a1bce48cbfb2e01559}{\texttt{31cc934}}. The links below identify the exact data and programs used for this paper. Download that revision and run the commands from its root directory.

\begin{description}[leftmargin=2.2em,style=nextline,itemsep=2pt,parsep=0pt]
\item[S6] \repodir{verification/n6}{verification/n6}: the sextic resultant, triangular factors, conjugacy identities, rational root intervals and recorded Magma input/output.
\item[S7] \repodir{verification/n7}{verification/n7}: ten rational boxes, contraction checks and twelve objective comparisons. The classification remains the external result \cite{fs}.
\item[S8] \repodir{verification/n8}{verification/n8}: exact matrix inputs in \texttt{data/global\_pair\_bound\_exact\_contacts.json}, algebraic and local tests, and the full-domain checker.
\item[S9] \repodir{verification/n9}{verification/n9}: \texttt{global\_pair\_bound\_exact\_contacts.json}, radial-moment bounds, boundary-contact tests and the full-domain checker.
\item[S10] \repofile{verification/n10/verify_counterexample.py}{verification/n10/verify\_counterexample.py}: the historical coordinates and the integer comparison in \eqref{historical}. The same directory also contains the alternative small-integer witness.
\item[G] \repofile{verification/ratio/check_examples.py}{verification/ratio/check\_examples.py}: the three integer sets in Table~\ref{tab:witnesses} and their strict distance-ratio comparisons.
\end{description}
The inventory \repofile{verification/files.sha256.json}{verification/files.sha256.json} binds the distributed files. First run \texttt{python verify.py} to audit it. Install \texttt{requirements.txt} with 64-bit Python 3.12 and leave assertions enabled. Then use

\begin{verbatim}
python verify.py --level check --case n6
python verify.py --level check --case n7
python verify.py --level replay --case n8
python verify.py --level replay --case n9
python verify.py --level check --case n10
python verify.py --level check --case ratio
\end{verbatim}
Each command must exit successfully. With \texttt{--json}, success is reported as \texttt{"passed": true}. The default audit checks stored evidence; the \texttt{replay} commands for eight and nine points rebuild the polynomials and start from the root region. Success requires positive matrix pivots, exact contact identities, all local tests, and zero unresolved boxes. The terminal counts must agree with Table~\ref{tab:coverage}. The terminal counts alone do not establish the pair inequality.

The six-point check reconstructs the resultant and rechecks the finite comparisons, while using the recorded elimination and dimension computation. The seven-point check compares twelve genuine classes and uses the published completeness theorem. These dependencies are distinct from calculations repeated by the programs. The \repofile{verification/README.md}{verification guide} explains the checks and their acceptance criteria.

\section{Integer witnesses}\label{integers}
The coordinates in Table~\ref{tab:witnesses} were chosen using the configurations illustrated in Friedman's catalogue as a guide \cite{friedman}. They are independent feasible examples. The last column lists the minimum and maximum squared distances obtained by enumerating all pairs. Repository material G checks these integers directly.

\begin{table}[H]
\centering
\caption{Integer configurations for the distance-ratio comparisons and their extreme squared distances.}\label{tab:witnesses}
\small\begin{tabularx}{\linewidth}{rXr}\toprule
$n$ & Coordinates & Min., max.\\\midrule
9 & $(15,28),$ $(5,22),$ $(25,22),$ $(15,16),$ $(0,11),$ $(30,11),$ $(15,4),$ $(4,0),$ $(26,0)$ & $136,925$\\
10 & $(14,27),$ $(5,22),$ $(22,22),$ $(14,17),$ $(0,14),$ $(27,14),$ $(11,8),$ $(2,4),$ $(25,4),$ $(16,0)$ & $89,733$\\
11 & $(9,17),$ $(15,16),$ $(4,14),$ $(11,12),$ $(6,9),$ $(16,9),$ $(0,8),$ $(11,6),$ $(4,4),$ $(16,3),$ $(9,1)$ & $29,289$\\\bottomrule
\end{tabularx}
\end{table}

\subsection{A comparison for all larger cardinalities}
For $m=n-2\ge10$, start with a regular $m$-gon of side one, and insert two points at distance $1/2$ on opposite sides of its center on a diameter. Its circumradius $r=1/(2\sin(\pi/m))>3/2$, and its inradius exceeds $1/2$. Vertex-to-inserted-point distances exceed one, so the separation is one and the diameter is the outer diameter $D_m$. The corresponding diameter of $S_n$ at separation one is $D_{m+1}$. For even $j$, $D_j=\csc(\pi/j)$; for odd $j$, $D_j=1/[2\sin(\pi/(2j))]$. In particular

\[
 D_m\le\csc(\pi/m)<(m+1)/\pi<D_{m+1}.
\]
The middle inequality follows from $\sin x\ge x-x^3/6$ and $\pi^2(m+1)<6m^2$ for $m\ge10$. This proves the comparison for all remaining cardinalities.

\subsection{The ten-point product check}
Repository material S10 stores the ten integer coordinates in Section~\ref{ten}, computes $q$, enumerates every unordered pair, and checks \eqref{historical} as an integer comparison. Its scope is solely the strict inequality $M_{10}>C_{10}$. Repository material G independently checks the minimum and maximum squared distances in Table~\ref{tab:witnesses}. Neither comparison requires a numerical optimizer.

\end{document}